\documentclass{aims}
\usepackage{amsmath}
  \usepackage{paralist}
  \usepackage{graphics} 
  \usepackage{epsfig} 
\usepackage{graphicx}  \usepackage{epstopdf}
 \usepackage[colorlinks=true]{hyperref}
\hypersetup{urlcolor=blue, citecolor=red}

\def\currentvolume{X}
 \def\currentissue{X}
  \def\currentyear{200X}
   \def\currentmonth{XX}
    \def\ppages{X--XX}

 \def\ocirc#1{\ifmmode\setbox0=\hbox{$#1$}\dimen0=\ht0
    \advance\dimen0 by1pt\rlap{\hbox to\wd0{\hss\raise\dimen0
    \hbox{\hskip.2em$\scriptscriptstyle\circ$}\hss}}#1\else
    {\accent"17 #1}\fi}
\def\R{{\mathbb R}}
\def\N{{\mathbb N}}

\newtheorem{theorem}{Theorem}[section]

\newtheorem{lemma}[theorem]{Lemma}

\theoremstyle{definition}

\newtheorem{remark}{Remark}

\title[Trace Theorem for H\" older domains]
      {A note on the Trace Theorem for domains which are locally subgraph of a H\" older continuous function}

\author[Boris Muha]{}


 \email{borism@math.hr}

\thanks{The author acknowledges post-doctoral support provided by the Texas Higher Education Coordinating Board, 
Advanced Research Program (ARP) grant number 003652-0023-2009 and MZOS grant number 0037-0693014-2765.}

\begin{document}
\maketitle

\centerline{\scshape Boris Muha }
\medskip
{\footnotesize
 \centerline{Department of Mathematics, Faculty of Natural Science, University of Zagreb,}
   \centerline{Bijeni\v cka 30, 10 000 Zagreb, Croatia}
} 

\bigskip

 \centerline{(Communicated by the associate editor name)}

\begin{abstract}
The purpose of this note is to prove a version of the Trace Theorem for domains which are locally subgraph of a H\" older continuous
function. More precisely, let $\eta\in C^{0,\alpha}(\omega)$, $0<\alpha<1$ and let $\Omega_{\eta}$ be a domain which is locally subgraph
of a function $\eta$. We prove that mapping $\gamma_{\eta}:u\mapsto u({\bf x},\eta({\bf x}))$ can be extended by continuity to a linear, continuous
mapping from $H^1(\Omega_{\eta})$ to $H^s(\omega)$, $s<\alpha/2$. This study is motivated by analysis of fluid-structure interaction problems.
\end{abstract}

\section{Introduction}
The Trace Theorem for Sobolev spaces is well-known and widely used in analysis of boundary and initial-boundary value problems in
partial differential equations. Usually, for the Trace Theorem to hold, the minimal assumption is that the domain
has a Lipshitz boundary (see e. g. \cite{ADA,Ding,Grisvard2}).
However, when studying weak solutions to a moving boundary fluid-structure interaction (FSI) problem, domains are not necessary
Lipshitz (see \cite{CDEM,CG,LenRuz,SunBorHyp2012,SunBorMulti}). 
FSI problems have many important applications (for example in biomechanics and aero-elasticity) 
and therefore have been extensively studied from the analytical, as well as numerical point of view, since the late 1990s (see e.g. 
\cite{CDEM,ChenShkoller,CG,Kuk,LenRuz,Leq13,BorSun} and the references within).
In FSI problems the fluid domain is unknown, given by an elastic deformation $\eta$, and therefore one cannot assume a priori any smoothness of the domain.
In \cite{CDEM,CG,LenRuz} an energy inequality implies $\eta\in H^2(\omega)$, 
$\omega\subset\R^2$. From the Sobolev embeddings one can see that in this case $\eta\in C^{0,\alpha}(\omega)$, $\alpha<1$, but
$\eta$ is not necessarily Lipschitz. Nevertheless, in Section 1.3 in \cite{CDEM}, and Section 1.3. in \cite{CG}, a version of the Trace Theorem
for such domains was proved, which enables the analysis of the considered FSI problems (see also \cite{LenRuz}, Section 2).

The proof of a version of the Trace Theorem in \cite{CG} (Lemma 2) relies on Sobolev embeddings theorems and the fact that $\eta\in H^2(\omega)$ and 
$\omega\subset\R^2$. 
Even though the techniques from \cite{CG} can be generalized to a broader class of Sobolev class boundaries, 
the result and techniques from \cite{CG} cannot be applied to some other cases of interest in FSI problems,
for example to the coupling of $2D$ fluid flow with the $1D$ wave equation, where we only have
$\eta\in H^1(\omega)$ (see \cite{SunBorHyp2012,SunBorMulti}) 
The purpose of this note is to fill that gap and generalize
that result for $\omega\subset\R^{n-1}$, $n>1$, and arbitrary H\" older continuous functions $\eta$. Hence,   
we prove a version of the Trace Theorem for a domain which is locally a subgraph of a H\" older continuous function.
We use real interpolation theory (see \cite{LionsMagenes}) and intrinsic norms for $H^s$ spaces, where $s$ in
not an integer.

\section{Notation and Preliminaries}

Let $n\in\N$, $n\geq 2$. Let $\omega\subset \R^{n-1}$ be a Lipschitz domain and let $0<\alpha<1$.
Furthermore, let $\eta$ satisfy the following conditions:
\begin{equation}\label{etacond}
\eta\in C^{0,\alpha}(\omega),\; \eta({\bf x})\geq \eta_{min}>0,\; {\bf x}\in \overline{\omega},\;  \eta_{|\partial\omega}=1.
\end{equation}
We consider the following domain
$$
\Omega_{\eta}=\{({\bf x},x_n):{\bf x}\in\omega,\; 0<x_n<\eta(x)\},
$$
with its upper boundary
$$
\Gamma_{\eta}=\{({\bf x},x_n):{\bf x}\in\omega,\; x_n=\eta(x)\}.
$$
We define the trace operator $\gamma_{\eta}:C(\overline{\Omega_{\eta}})\to C(\omega)$
\begin{equation}\label{trace}
(\gamma_{\eta}u)({\bf x})=u({\bf x},\eta({\bf x})),\quad {\bf x}\in \omega,\;u\in C^0(\overline{\Omega_{\eta}}).
\end{equation}
In \cite{CDEM} (Lemma 1) it has been proven that $\gamma_{\eta}$ can be extended by continuity to 
an operator $\gamma_{\eta}:H^1(\Omega_{\eta})\to L^2(\omega)$.
This result holds with an assumption that $\eta$ is only continuous.
Our goal is to extend this result in a way to show that Im$(\gamma_{\eta})$ is a subspace of $H^s(\omega)$, for some $s>0$,  when $\eta$ is 
a H\" older continuous function.
\begin{remark}
Notice that $\gamma_{\eta}$ is not a classical trace operator because $\gamma_{\eta}(u)$ is a function defined on $\omega$, whereas the classical trace
would be defined on the upper part of the boundary, $\Gamma_{\eta}$. However, this version of a trace operator is exactly what one needs in
analysis of FSI problems. Namely, in the FSI setting the Trace Theorem is applied to fluid velocity which, at the interface,
equals the structure velocity, where the structure velocity 
is defined on a Lagrangian domain (in our notation $\omega$).
\end{remark}

The Sobolev space $H^s(\omega)$, $0<s<1$ is defined by the real interpolation method (see \cite{ADA,LionsMagenes}).
However, $H^s(\omega)$ can be equipped with an equivalent, intrinsic norm (see for example \cite{ADA,Grisvard2}) which is also used in \cite{Ding}
\begin{equation}\label{Hsnormdef}
\displaystyle{\|u\|^2_{H^s(\omega)}=\|u\|^2_{L^2(\omega)}+\int_{\omega\times\omega}\frac{|u({\bf x}_1)-u({\bf x}_2)|^2}{|{\bf x}_1-{\bf x}_2|^{n-1+2s}}d{\bf x}_1d{\bf x}_2,}
\end{equation}
where $0<s<1$.
\section{Statement and Proof of the result}

\begin{theorem}\label{main}
Let $\alpha<1$ and let $\eta$ be such that conditions (\ref{etacond}) are satisfied.
Then operator $\gamma_{\eta}$, defined by (\ref{trace}), can be extended by continuity to a linear operator from $H^1(\Omega_{\eta})$ to 
$H^s(\omega)$, $0\leq s< \frac{\alpha}{2}$.
\end{theorem}
\proof
We split the main part of the proof into two Lemmas. The main idea of the proof is to transform a function defined on $\Omega_{\eta}$ to a
function defined on $\omega\times (0,1)$ and to apply classical Trace Theorem to a function defined on the domain $\omega\times (0,1)$. 
Throughout this proof $C$ will denote a generic positive constant that depends only on $\omega$, $\eta$ and $\alpha$.

Let $u\in H^1(\Omega_{\eta})$. Define 
\begin{equation}\label{trans}
\bar{u}({\bf x},t)=u({\bf x},\eta({\bf x})t),\; {\bf x}\in\omega,\; t\in [0,1].
\end{equation}
Let us define function space (see \cite{LionsMagenes}, p. 10):
$$
W(0,1;s)=\{f:f\in L^2(0,1;H^s(\omega)),\; \partial_t f\in L^2(0,1;L^2(\omega))\},
$$
where $0<s<1$. Our goal is to prove $\bar{u}\in W(0,1;s)$. However, before that we need to prove the following
technical Lemma:
\begin{lemma}\label{krivulja}
For every ${\bf x}_0$, ${\bf x}_1\in\omega$, there exists a piece-wise smooth curve parameterized by
$$\Theta_{x_0,x_1}:[0,2]\to\Omega_{\eta}$$
such that $\Theta_{x_0,x_1}(0)=({\bf x}_0,\eta({\bf x}_0))$, $\Theta_{x_0,x_1}(2)=({\bf x}_1,\eta({\bf x}_1))$ and
\begin{equation}\label{CurveDer}
|\Theta_{x_0,x_1}'(r)|\leq C|{\bf x}_1-{\bf x}_0|^{\alpha},\quad {\rm a.\;e.}\; r\in [0,2],
\end{equation}
where $C$ does not depend on ${\bf x}_0$, ${\bf x}_1$.
\end{lemma}
\proof
First we define ${\bf x}_r$ as a convex combination of ${\bf x}_0$ and ${\bf x}_1$:
$$
{\bf x}_r=(1-r^{1/\alpha}){\bf x}_0+r^{1/\alpha}{\bf x}_1={\bf x}_0+r^{1/\alpha}({\bf x}_1-{\bf x}_0),\quad r\in [0,1].
$$
Furthermore we define ${\bf y}_r$ in the following way:
$$
y_r=\eta({\bf x}_0)-\|\eta\|_{C^{0,\alpha}(\omega)}|{\bf x}_r-{\bf x}_0|^{\alpha}=
\eta({\bf x}_0)-\|\eta\|_{C^{0,\alpha}(\omega)}r|{\bf x}_1-{\bf x}_0|^{\alpha},\; r\in[0,1].
$$
By using H\" older continuity of $\eta$ we get
\begin{equation}\label{yr}
y_r\leq\eta({\bf x}_r),\quad r\in [0,1].
\end{equation}
Therefore curve $({\bf x}_r,y_r)$ stays bellow the graph of $\eta$ for $r\in [0,1]$. Now, let us consider whether
this curve intersects the hyper-plane $x_n=\eta_{min}$. Since $y_r$ is a strictly decreasing function in $r$, 
we distinguish between  the two separate cases.
\vskip 0.1in 
\noindent
{\bf Case 1}: $y_r\geq \eta_{min}$, $r\in [0,1]$.
We define $\Theta_{x_0,x_1}$ in the following way:
\begin{equation}\label{case1}
\Theta_{x_0,x_1}(r)=\left \{\begin{array}{lcr} ({\bf x}_r,y_r) &,& 0\leq r\leq 1, \\ \\
({\bf x}_1,(2-r)y_1+(r-1)\eta({\bf x}_1)) &,& 1<r\leq 2. \end{array} \right .
\end{equation}
From (\ref{yr}), the definition of $\Theta_{x_0,x_1}$ (\ref{case1}) and the definition of $\Omega_{\eta}$ it follows immediately that
$\Theta_{x_0,x_1}(0)=({\bf x}_0,\eta({\bf x}_0))$, $\Theta_{x_1,x_2}(2)=({\bf x}_1,\eta({\bf x}_1))$ and
$\Theta_{x_0,x_1}(r)\in\Omega_{\eta},\; r\in[0,2]$. Therefore it only remains to prove (\ref{CurveDer}). We calculate
$$
\Theta_{x_0,x_1}'(r)=\left \{\begin{array}{lcr} 
(\frac{1}{\alpha}r^{1/\alpha-1}({\bf x}_1-{\bf x}_0),-\|\eta\|_{C^{0,\alpha}(\omega)}|{\bf x}_1-{\bf x}_0|^{\alpha}) &,& 0\leq r\leq 1, \\ \\
(0,\eta({\bf x}_1)-y_1) &,& 1<r\leq 2. 
\end{array} 
\right .
$$
Since $\omega$ is bounded, we can take $C\geq\|\eta\|_{C^{0,\alpha}(\omega)}$ such that
$$
|{\bf x}-{\bf y}|\leq C |{\bf x}-{\bf y}|^{\alpha},\quad {\bf x},\; {\bf y}\in\omega.
$$
Using this observation we can get an estimate:
$$
|\Theta_{x_0,x_1}'(r)|\leq C |{\bf x}_0-{\bf x}_1|^{\alpha},\quad r\in [0,1).
$$
Furthermore, analogously using the definition of $y_r$ and $\eta\in C^{0,\alpha}(\omega)$ we have
$$
|\eta({\bf x}_1)-y_1|\leq |\eta({\bf x}_1)-\eta({\bf x}_0)|+\|\eta\|_{C^{0,\alpha}(\omega)}r|{\bf x}_1-{\bf x}_0|^{\alpha}
\leq C |{\bf x}_0-{\bf x}_1|^{\alpha}.
$$
Therefore,  (\ref{CurveDer}) is proven.
\vskip 0.1in 
\noindent
{\bf Case 2}: There exists $r_0\in (0,1)$ such that $y_r=\eta_{min}$. In this case we define
$\Theta_{x_0,x_1}$ in the following way:
\begin{equation}\label{case2}
\Theta_{x_0,x_1}(r)=\left \{\begin{array}{lcr} ({\bf x}_r,y_r) &,& 0\leq r\leq r_0, \\ \\
({\bf x}_r,\eta_{min}) &,& r_0< r\leq 1, \\ \\
({\bf x}_1,(2-r)\eta_{min}+(r-1)\eta({\bf x}_1)) &,& 1<r\leq 2. \end{array} \right .
\end{equation}
Analogous calculation as in Case 1 shows that estimate (\ref{CurveDer}) is valid in this case as well.
This completes the proof of the Lemma.
\qed

Now we are ready to prove the following lemma:
\begin{lemma}\label{translemma}
Let $u\in H^1(\Omega_{\eta})$ and let $0<s<\alpha$. Then ${\bar u}\in W(0,1;s)$, where ${\bar u}$ is defined by formula (\ref{trans}).
\end{lemma}
\proof
Let us first take $u\in C^{\infty}_c(\R^n)$. For ${\bf x}_1, {\bf x}_2\in\omega$, $t\in (0,1)$ we have
$$
|{\bar u}({\bf x}_1,t)-{\bar u}({\bf x}_2,t)|=|u({\bf x}_1,\eta({\bf x}_1)t)-u({\bf x}_2,\eta({\bf x}_2)t)|
$$
Notice that $t\eta\in C^{0,\alpha}(\omega)$ and therefore we can apply Lemma \ref{krivulja} 
to function $t\eta$ (we just need to replace $\eta_{min}$ with $t\eta_{min}$ in the proof of the Lemma \ref{krivulja})
to get $\Phi_{x_1,x_2}^t:[0,2]\to\Omega_{\eta}$ such that:
$$
\begin{array}{c}
\displaystyle{\Theta^t_{x_1,x_2}(0)=({\bf x}_1,\eta({\bf x}_1)t),\quad \Theta^t_{x_1,x_2}(2)=({\bf x}_2,\eta({\bf x}_2)t)},
\\ \\
|\displaystyle{\frac{d}{dr}\Theta^t_{x_1,x_2}(r)|\leq C|{\bf x}_1-{\bf x}_2|^{\alpha},\quad {\rm a.\;e.}\; r\in [0,2],}
\end{array}
$$
where $C$ does not depend on ${\bf x}_1$, ${\bf x}_2$ and $t$.
Define
$$
f^t_{x_1,x_2}(r)=u(\Theta^t_{x_1,x_2}(r)),\; r\in [0,2].
$$
Now we have
\begin{equation}\label{Est1}
\begin{array}{c}
\displaystyle{|u({\bf x}_1,\eta({\bf x}_1)t)-u({\bf x}_2,\eta({\bf x}_2)t)|^2=|\int_0^2\frac{d}{dr}f^t_{x_1,x_2}(r)dr|^2}
\\ \\
\displaystyle{\leq \|\frac{d}{dr}\Theta^t_{x_1,x_2}(r)\|_{L^{\infty}(0,2)}^2\int_0^2|\nabla u(\Theta^t_{x_1,x_2}(r))|^2dr
}
\\ \\
\displaystyle{\leq C|{\bf x}_1-{\bf x}_2|^{2\alpha}\int_0^2|\nabla u(\Theta^t_{x_1,x_2}(r))|^2dr.}
\end{array}
\end{equation}
Using (\ref{Est1}) we get the following estimates:
\begin{equation}\label{Hsnorm}
\begin{array}{c}
\displaystyle{\|{\bar u}\|_{L^2(0,1:H^s(\omega))}^2=\int_0^1\|\bar{u}(.,t)\|^2_{H^s(\omega)}dt
=\int_0^1dt\int_{\omega\times\omega}\frac{|\bar{u}({\bf x}_1,t)-\bar{u}({\bf x}_2,t)|^2}{|{\bf x}_1-{\bf x}_2|^{n-1+2s}}d{\bf x}_1d{\bf x}_2}
\\ \\
\displaystyle{\leq C \int_0^1dt\int_{\omega\times\omega}\frac{d{\bf x}_1d{\bf x}_2}{|{\bf x}_1-{\bf x}_2|^{n-1+2(s-\alpha)}}
\int_0^2|\nabla u(\Theta^t_{x_1,x_2}(r))|^2dr.}
\\ \\
\displaystyle{\leq C\|\nabla u\|^2_{L^2(\Omega_{\eta})}\int_{\omega\times\omega}\frac{d{\bf x}_1d{\bf x}_2}{|{\bf x}_1-{\bf x}_2|^{n-1+2(s-\alpha)}}}.
\end{array}
\end{equation}

To estimate the last integral in (\ref{Hsnorm}), we introduce a new variable ${\bf h}={\bf x}_1-{\bf x}_2$ and the change of variables
$({\bf x}_1,{\bf x}_2)\mapsto ({\bf h},{\bf x}_2)$ to get:
\begin{equation}\label{termB}
\displaystyle{\int_{\omega\times\omega}\frac{d{\bf x}_1d{\bf x}_2}{|{\bf x}_1-{\bf x}_2|^{n-1+2(s-\alpha)}}\leq C\int_{-R}^R\frac{dh}{|h|^{1+2(s-\alpha)}}},
\end{equation}
where $R={\rm diam}(\omega)$.
Recall that $s<\alpha<1$. Therefore by combining (\ref{Hsnorm}) and (\ref{termB}), we get:
\begin{equation}\label{HSglatke}
\|\bar{u}\|_{L^2(0,1:H^s(\omega))}\leq C\|u\|_{H^1(\Omega_{\eta})},\quad u\in C^{\infty}_c(\R^n)
\end{equation}
Since $C^{\infty}_c(\R^n)$ is dense in $H^1(\Omega_{\eta})$ (see \cite{ADA}, Thm 2, p. 54 with a slight modification near
$\partial\omega\times\{1\}$, see also \cite{CDEM}, proof of Lemma 1 and \cite{LenRuz}, Prop A.1.), by a density argument
we have 
$$
\bar{u}\in L^2(0,1;H^s(\omega)),\quad u\in H^1(\Omega_{\eta}).
$$
Now, it only remains to prove $\partial_t\bar{u}\in L^2((0,1)\times\omega)$. However, this can be proven with direct calculation by
using the chain rule:
$$
\partial_t\bar{u}({\bf x},t)=\eta({\bf x})\partial_{x_n} u({\bf x},\eta({\bf x})t).
$$
Since $\eta$ is H\" older continuous on $\omega$, from the above formula we have $\partial_t\bar{u}\in L^2((0,1)\times\omega)$ which
completes the proof of the Lemma.
\qed

\hskip 0.1in
Now we use continuity properties of $W(0,1;s)$ (\cite{LionsMagenes}, p. 19, Thm 3.1.), i.e.
$$
W(0,1;s)\hookrightarrow C([0,T];H^{s/2}(\omega)),
$$
where this injection is continuous.
Therefore, from Lemma \ref{translemma} we have
\begin{equation}\label{cont}
\bar{u}\in C([0,T];H^{s/2}(\omega)),\quad u\in H^1(\Omega_{\eta}).
\end{equation}
We finish the proof by noticing that $\gamma_{\eta}(u)=\bar{u}(.,1)$.
\qed

\begin{remark}
In \cite{CG}, Lemma 2, a special case of Theorem \ref{main} was proved. Namely, for $n=3$ and $\eta\in H^2(\omega)$ it was proved that 
$\gamma_{\eta}$ is a continuous operator from $H^1(\Omega_{\eta})$ to $H^s(\omega)$, $0\leq s<\frac 1 2$. This result follows from Theorem
\ref{main} because of the Sobolev imbedding $H^2(\omega)\hookrightarrow C^{0,\alpha}(\omega)$, $\alpha<1$. However, the techniques from \cite{CG}
rely on Sobolev embeddings and the fact that $\nabla\eta$ is more regular then $L^2(\omega)$ and therefore, cannot be extended for the case
of arbitrary H\" older continuous functions.
\end{remark}

\medskip
Received xxxx 20xx; revised xxxx 20xx.
\medskip

\end{document}